\documentclass[11pt]{amsart}
\usepackage{url}
\newtheorem{theorem}{\bf Theorem}
\newtheorem{definition}[theorem]{\bf Definition}
\newtheorem{remark}[theorem]{\bf Remark}

\newtheorem{proposition}[theorem]{\bf Proposition}
\newtheorem{corollary}[theorem]{\bf Corollary}
\begin{document}
\noindent
{\small Topology Atlas Invited Contributions \textbf{9} no.~2 (2004) 7 pp.}
\vskip 0.5in
\title{A survey of topological work at CEOL}
\author{Homeira Pajoohesh}
\address{Homeira Pajoohesh\\
Boole Centre for Research in Informatics}
\author{M. P. Schellekens}
\address{M. P. Schellekens\\
Centre for Efficiency Oriented Languages (CEOL)\\
Department of Computer Science\\
National University of Ireland, Cork}
\urladdr{http://www.ceol.ucc.ie/}
\thanks{CEOL is funded by Science Foundation Ireland, Investigator Award 
02/IN.1/181}
\begin{abstract} 
We present an overview of ongoing work at the Centre for
Efficiency-Oriented languages (CEOL), with a focus on topological aspects.
CEOL researchers are engaged in designing a new Real-Time Language to
improve software timing. The centre broadly focuses on bridging Semantics
and Complexity and unites researchers with expertise in Semantics of
Programming Languages, Real-Time Languages, Compiler Design and Graph
Based Algorithms. CEOL aims to narrow the gap between Worst Case Execution
Time analysis and Average Case Execution Time analysis for Real-Time
languages and its longer term goal is the development of ACETT, an Average
Case Execution Time Tool. This research work is of crucial interest to
industry, given that real-time software is widely used in a variety of
applications, such as chemical plants, satellite communications, the space
industry, telephone exchanges, medical equipment, the motor industry, etc.
Topological work at CEOL focuses on the exploration of Quantitative
Domains, semivaluations, partial metrics and their applications. We give
an overview of prior results obtained at CEOL in this area and of current
work on relating the notion of \emph{balance} of algorithms to running time
and of an exploration of semivaluations in relation to algorithmic running
time.
\end{abstract}
\maketitle

\section{Quantitative Domain Theory} 

Domain Theory, a formal basis for the semantics of programming languages,
originated in work by Dana Scott in the mid-1960s. Models for various
types of programming languages, including imperative, functional,
nondeterministic and probabilistic languages, have been extensively
studied. Quantitative Domain Theory forms a new branch of Domain Theory
and has undergone active research in the past three decades. The field
involves both the semantics of programming languages and Topology.

In order to reconcile two alternative approaches to Domain Theory, the
order-theoretic approach and the metric approach, Michael Smyth pioneered
at Imperial College the use of methods from the field of Non-Symmetric
Topology \cite{sm1}. This field traditionally studies quasi-metrics
which are obtained from classical metrics by removing the symmetry
requirement.  Hence, for quasi-metrics, the distance from a given point to
a second point need not be the same as the converse distance. A simple
example of a quasi-metric is the $0$--$1$ encoding of a partial order
which defines the distance between two points $x$ and $y$ to be $0$ in
case $x$ is below $y$ in the order and $1$ otherwise.

Recent developments in Domain Theory indicate that additional concepts are
required in order to develop the corresponding applications. We focus on
applications related to complexity analysis. For other examples we refer
the reader to \cite{sch8}. An extensive series of papers has been
published in this area, by the second listed author in collaboration with
Salvador Romaguera (e.g., \cite{rs2, rs, rs1, ms}). We will focus on the
connections between topology and complexity analysis below.

Each of these applications involve real distances in some sense, and hence
the adjective quantitative is used as opposed to the adjective
qualitative, indicating the traditional order-theoretic approach.

To develop Quantitative Domain Theory, we favour partial metrics
(pmetrics) whose category is isomorphic to the category of weightable
quasi-metric spaces, introduced by Steve Matthews at Warwick, \cite{M94}.  
Pmetrics are obtained from classical metrics by removing the requirement
that self distance be zero.  We recall that the qmetric $q\colon X\times
X\to \mathbb{R}^{+}_{0}$ is weightable if there is a function $w\colon
X\to \mathbb{R}^{+}_{0}$ such that for every $x,y \in X$,
$q(x,y)+w(x)=q(y,x)+w(y)$.  It is well known that this approach
considerably simplifies topological completions as remarked in
\cite{sch8}. To simplify matters even further, our approach is to view
Quantitative Domain Theory in first instance as an extension of
traditional Domain Theory via a minimal fundamental concept: that of a
semivaluation \cite{ms4, sch8}. A semivaluation is a novel
mathematical notion which generalises the fruitful notion of a valuation
on a lattice to the context of semilattices. It essentially emerged from a
study of the self-distance for partial metrics.

\begin{definition}
If $(X,\preceq)$ is a meet semilattice then a function 
$f\colon {X \to \mathbb{R}^{+}_{0}}$ 
is a {\it meet valuation} iff
$$\forall x,y,z\in X\ f(x \sqcap z) \geq f(x \sqcap y) + f(y \sqcap z) - f(y) $$ 
and $f$ is {\it meet co-valuation} iff 
$$\forall x,y,z\in X\ f(x \sqcap z) \leq f(x \sqcap y) + f(y \sqcap z) - f(y).$$
\end{definition}

\begin{definition}
If $(X,\preceq)$ is a join semilattice then a function 
$f\colon {X \to \mathbb{R}^{+}_{0}}$ 
is a {\it join  valuation} iff 
$$\forall x,y,z\in X\ f(x \sqcup z) \leq  f(x \sqcup y) + f(y \sqcup z) -f(y)$$ 
and $f$ is {\it join co-valuation} iff
$$\forall x,y,z\in X\ f(x \sqcup z) \geq f(x \sqcup y) + f(y \sqcup z) -f(y).$$
\end{definition}

\begin{definition}
A function is a {\it semivaluation} if it is either a join valuation or a 
meet valuation. A join (meet) valuation space is a join (meet)   
semilattice equipped with a join (meet) valuation. A {\it semivaluation 
space} is a semilattice equipped with a semivaluation.
\end{definition}

\break
\begin{proposition} Let $L$ be a lattice.
\begin{enumerate}
\item
A function $f\colon L \to \mathbb{R}^{+}_{0}$ is a join valuation if and 
only if it is increasing and satisfies join-modularity, i.e., 
$$f(x \sqcup z) + f(x \sqcap z) \leq f(x) + f(z).$$
\item
A function $f\colon L \to \mathbb{R}^{+}_{0}$ is a meet valuation if and 
only if it is increasing and satisfies meet-modularity, i.e., 
$$f(x \sqcup z) + f(x \sqcap z) \geq f(x) + f(z).$$
\end{enumerate}
\end{proposition}

\begin{corollary} 
A function on a lattice is a valuation iff it is a join valuation and a 
meet valuation. A function on a lattice is a co-valuation iff  it is a 
join co-valuation and a meet co-valuation.
\end{corollary}

\begin{proposition} 
If a function $f\colon L \to \mathbb{R}^{+}_{0}$ is join co-valuation 
then $p \colon L\times L \to \mathbb{R}^{+}_{0}$ with  
$p(x,y)=f(x)-f(x\vee y)$ is a pmetric. 
\end{proposition}

Proposition 4 clearly motivates the fact that semivaluations provide a
natural generalisation of valuations from the context of lattices to the
context of semilattices. We refer the reader to \cite{ms4} for the
correspondence theorems which link partial metrics to semivaluations.

The semivaluation approach has the advantage that it allows for a uniform
presentation of the traditional quantitative domain theoretic structures
and applications, as for instance the totally bounded Scott domains of
Smyth and the partial metric spaces of Matthews.

A topological problem stated by the Hans-Peter K\"unzi essentially
required the mathematical characterization of partial metrics. Such a
characterization has been obtained based on the notion of a semivaluation
\cite{sch8}. This notion was directly motivated by Computer Science
examples. Hence the result forms an example of recent developments where
the mathematical area of Topology is influenced by Computer Science.
Traditionally the influence has largely been in the opposite direction.

The benefit to Computer Science is that semivaluations allow for the
introduction of a suitable notion of a \emph{quantitative domain} which 
can serve to develop models for the above mentioned applications.

Below we indicate some results currently under investigation at CEOL in
relation to partial metrics and semivaluations on binary trees, complexity
analysis and on the imbalance lattice of binary trees.

\section{Binary trees, semivaluations and complexity analysis} 

\subsection{A complexity pmetric on decision trees}

\begin{definition}
Binary trees are downwards-directed finite trees with a distinguished root
node, in which every non-root node has either exactly two children or no
children (i.e., it is a leaf). The placement of the nodes at each level is
not significant, so a tree is determined (up to permutation) by the
path-lengths of its leaves (i.e. the lengths of the paths (=branches) from
the root node to the leaves). Thus we can represent equivalence classes of
the rooted binary trees with $n$ leaves by sequences of $n$ non-negative
integers, which give the path-length of each leaf. For example, the
path-length sequence $\langle 1 3 3 4 4 4 4 \rangle$ represents a binary
tree with $n=7$ leaves, of which one has path-length $1$, two have
path-length $3$, and four have path-length $4$. This is shown in Figure 1.
\end{definition}

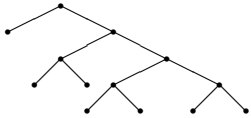
\begin{figure}
\begin{picture}(90,40)(0,0)
\put(20,40){\circle*{2}}
\put(20,40){\line(-2,-1){20}}
\put(20,40){\line(2,-1){20}}
\put(0,30){\circle*{2}}
\put(40,30){\circle*{2}}
\put(40,30){\line(-2,-1){20}}
\put(40,30){\line(2,-1){20}}
\put(20,20){\circle*{2}}
\put(20,20){\line(1,-1){10}}
\put(20,20){\line(-1,-1){10}}
\put(60,20){\circle*{2}}
\put(60,20){\line(-2,-1){20}}
\put(60,20){\line(2,-1){20}}
\put(30,10){\circle*{2}}
\put(10,10){\circle*{2}}
\put(40,10){\circle*{2}}
\put(40,10){\line(1,-1){10}}
\put(40,10){\line(-1,-1){10}}
\put(80,10){\circle*{2}}
\put(80,10){\line(1,-1){10}}
\put(80,10){\line(-1,-1){10}}
\put(90,0){\circle*{2}}
\put(70,0){\circle*{2}}
\put(50,0){\circle*{2}}
\put(30,0){\circle*{2}}
\end{picture}
\caption{The binary tree $\langle 1 3 3 4 4 4 4 \rangle$}
\end{figure}

We can interpret a binary tree as a $\vee$-semilattice and conversely
define a binary tree as a semilattice with some additional properties.
This leads to a characterization of binary trees in terms of semilattices.  
Our interest in binary trees stems from the fact that in order to carry
out the running time analysis of so-called \emph{comparison-based
algorithms}, i.e., algorithms in which every action is ultimately based
on a prior comparison between two elements, the notion of a decision tree
is a fundamental tool \cite{ah}. Typical examples are so-called
\emph{Sorting Algorithms} which sort a given list of numbers in increasing
order. Decision trees are binary trees representing the comparisons
carried out during the computation of a comparison-based algorithm. The
distance from the root (\emph{input}) to a leaf (\emph{output}) gives 
the comparison time for the algorithm (i.e., the total number of comparisons)
to compute the output corresponding to the given leaf.

\begin{proposition} 
The distance from any node of a binary tree to the root is a join  
co-valuation. 
\end{proposition}

\begin{proof} 
We sketch the proof. We proceed by induction on the number of leaves. 
It suffices to show that if $T\in T_n$, where $T_{n}$ is the set of binary 
trees with $n$ leaves, and $x,y,z\in T$ then
\begin{equation}
l(x\vee z)\geq l(x\vee y)+l(y\vee z)-l(y).\tag{$\star$}
\end{equation}

Here $l(x)$, the level of $x$, is the distance from $x$ to the root. 

We show that ($\star$) holds for leaves. The result for general nodes then
follows.  This is true for $n=2$. We assume that this is true for every
$k$ less than $n$ and we have to show it still holds for $n$.  Notice that
if the level of none of them is the maximum level then we can remove this
maximum level and by induction ($\star$) is true for $x,y,z$.

Now if we have some leaves at levels other than the levels of $x,y$ and
$z$, we can remove the leaves and use induction. Thus the problem is 
reduced to the case where there are no other levels except $l(x),l(y)$ 
and $l(z)$.

On the other hand notice, that if $x\vee y\vee z$ is not the root then 
again we can reduce the problem by induction and prove that ($\star$) is 
true. So we assume that this join is the root. We have to investigate 
three cases: 
\begin{enumerate}
\item
Three of them are at the same level.
\item
Two of them are at the same level and the last one is on a different level. 
\item
They are on three pairwise different levels. 
\end{enumerate}

We prove ($\star$) in the first case and for the other cases the proof is
similar. 

Since three of them are at the same level they have to be in the last
level. Notice that the tree divides via two branches from the root in two
subtrees. Since the join of the three elements is the root, three of them
can't be in the same half, so two of them are in the same half and the
other is the other half. Now we have to investigate two cases, $x,y$ are
in the same half and $z$ is in the other half and the second case where
$x,z$ are in the same half and $y$ is in another half.

If $x,y$ are in the same half and $z$ is in the other half then both
$x\vee z$ and $y\vee z$ are the root and so $l(x\vee z)=l(y\vee z)=0$ and
since $l$ is order reversing $l(x\vee y)\geq l(y)$ and ($\star$) trivially
holds.

Also if $x,z$ be in the same half and $y$ in another half then the 
argument is similar and the proof is complete. 
\end{proof}

\begin{remark}
It is easy to verify that the proof extends to the case of trees for 
which every node has at most two children.
\end{remark}

By Proposition 6 we have a partial metric (pmetric) on each tree
\cite{sch8}. When we consider this pmetric in the context of decision
trees, we conclude that the maximum of the weights of the leaves gives the
worst-case running-time, while the average of the weights of the leaves
gives the average comparison time.

Every partial metric on a set induces two topologies on the set, see
\cite{KMP03}.  Here the base of one of these two topologies is the set of
upper sets and the base of the other one is the set of down sets. But
looking at the down sets is interesting, because by the nature of binary
trees, every down set is a binary tree with less leaves. Then the down set
$T$ is called chain-open if the set of down sets which are subsets of $T$
is a chain with respect to inclusion. We can show that the set of maximal
chain open sets is a base for a binary tree and we can characterize the
binary tree of insertion sort.

\subsection{On Balance and Algorithmic Running Time}

In \cite{ah} it is argued that Divide \& Conquer techniques that make a
\emph{balanced} division in general lead to faster algorithms. We use 
the so called imbalance lattice in the following to put this intuition on 
a formal basis.

\begin{definition}
If $x = \langle x_1,\dots,x_n\rangle, y = \langle y_1,\dots,y_n\rangle \in 
T_n$ then we say that $x$ is more balanced than $y$ if 
$\sum_{i=1}^{n} x_i\leq \sum_{i=1}^{n} y_i$. 
\end{definition}

\begin{definition} 
If $x,y\in T_n$ then we define $x\leq y$ if there are 
$l_1,\dots,l_m \in T_n$ such that $y=l_1,x=l_m$ and for each $i$, 
$1\leq i<m$, $l_{i+1}$ is more balanced than $l_i$. 
\end{definition}

\begin{proposition}[\cite{pr}]
$T_{n}$ with the above order is a lattice.
\end{proposition} 

For any $n$ we define a pmetric on $T_n$ by regarding the order on it via
the classical technique of \cite{kv} involving shortest paths. Then for
the two topologies that we mentioned in the last paragraph, one of them is
the set of upper sets and the other is the set of lower sets on it and the
Lawson topology, which is the join of these topologies, is the discrete
topology \cite{ps}.

Finally we mention that we investigated the use of the imbalance lattice
\cite{pr} in the context of path-length sequences for decision trees.
This forms a novel exploration of this lattice which priorly was solely
restricted to the study of balance of trees, independently from running
time. In this way we obtained a link between the balance of the decision
tree of an algorithm and the speed of the algorithm.

Recent CEOL work \cite{ops} has shown that the balance lattice is a useful
tool to illustrate that \emph{more balanced algorithms} are faster. This
intuition was raised in the literature \cite{ah}. We illustrated in
particular that Mergesort can be interpreted as a more balanced version of
Insertionsort, leading to an alternative proof that the first algorithm is
faster than the second. This corresponds to \cite{ah} where it is argued
that Divide \& Conquer algorithms that make a \emph{balanced} division in
general are faster algorithms \cite{ops}. Our work has put this intuition
on a formal basis.

\providecommand{\bysame}{\leavevmode\hbox to3em{\hrulefill}\thinspace}
\providecommand{\MR}{\relax\ifhmode\unskip\space\fi MR }

\end{document}